\documentclass{article}
\pdfoutput=1

\usepackage[a4paper,vmargin=30mm,hmargin=20mm]{geometry}
\usepackage[dvipsnames]{xcolor}
\usepackage[T1]{fontenc}
\usepackage{microtype}
\usepackage{amsmath}
\usepackage{amssymb}
\usepackage{booktabs}
\usepackage{natbib}
\setcitestyle{round}

\usepackage[unicode,bookmarks,colorlinks,citecolor=blue,urlcolor=BrickRed]{hyperref}
\hypersetup{pdftitle={Analytical computation of moderate-degree fully-symmetric cubature rules on the triangle},pdfauthor={Stefanos-Aldo Papanicolopulos},pdfkeywords={cubature, triangle, fully symmetric rules, symmetric polynomials},pdfsubject={Cubature}}

\newcommand{\doi}[1]{doi: \href{http://dx.doi.org/#1}{#1}}

\begin{document}

\title{Analytical computation of moderate-degree\\ fully-symmetric cubature rules on the triangle}

\author{Stefanos-Aldo~Papanicolopulos%
\thanks{Department of Mechanics, National Technical University of Athens, Zografou 15773, Greece, email: stefanos@mechan.ntua.gr}}

\date{November 16, 2011}

\maketitle

\begin{abstract}
A method is developed to compute analytically fully symmetric cubature rules on the triangle by using symmetric polynomials to express the two kinds of invariance inherent in these rules. Rules of degree up to 15, some of them new and of good quality, are computed and presented.
\end{abstract}


\section{Introduction}
  \label{sec:intro}

The term ``cubature'' indicates the numerical computation of a multiple integral. This is an important topic in many different disciplines, with a correspondingly large body of literature. A description of the different kinds of cubature rules that exist, as well as of the mathematics used to derive them, is given in the classical book of \citet{Stroud1971}, with more updated information to be found, among others, in \citep{Cools1997} and in chapter~6 of~\citep{KrommerUeberhuber1998}. \Citet{Stroud1971} also presents a compilation of known (at the time) cubature rules, while newer rules are catalogued in \citep{CoolsRabinowitz1993,Cools1999} and online at the Encyclopedia of Cubature Formulas~\citep{Cools2003}.

A commonly used method to derive specific cubature rules is based on moment equations and invariant theory \citep[pp.~170--182]{KrommerUeberhuber1998}. This method, which will be used in the present paper, exploits symmetries and invariant theory to set up a non-linear system of equations, whose unknowns are the positions and weights of the integration points. The use of invariants, together with appropriate analytical computations, can lead to a significant simplification of the system of equations, which however in most cases still has to be solved numerically.

Although appropriate iterative numerical methods have been successfully used to obtain \emph{individual} numerical solutions to the aforementioned system of equations, obtaining a solution in this way provides no information on its uniqueness. Conversely, inability to obtain a solution does not prove its inexistence (though it is a strong indication, when sufficiently robust numerical methods are employed). It is thus interesting and useful to be able to compute analytically all the solutions for a cubature rule.

In this paper we focus on fully symmetric cubature rules on the triangle and provide results of the analytical computation for cubature rules of moderate degree, extending significantly the analytical results given by \citet{LynessJespersen1975}. The theory of symmetric polynomials~\citep{Macdonald1998} is used in the generation of the non-linear system of equations, resulting in a straightforward way to formulate a system which is simpler than those resulting from the use of polar coordinates~\citep{LynessJespersen1975,BerntsenEspelid1990} or Cartesian coordinates~\citep{WandzuraXiao2003}. Additionally, symmetric polynomials are also used to significantly simplify the computation and presentation of the analytical solution.

\section{Symmetric polynomials}

Although the formulation presented here is actually based on invariant theory, the relevant theory is not used directly but is ``implied'' by using the theory of symmetric polynomials~\citep{Macdonald1998}. As we will see in the following, the use of symmetric polynomials provides a concise formulation of the non-linear system of equations, while also leading to simpler computation and presentation of the solution.

A \emph{symmetric polynomial} is a multivariate polynomial in $n$ variables, say $x_1,x_2,\ldots,x_n$, which is invariant under any permutation of its variables. We define the \emph{elementary symmetric polynomials} $\tilde{x}_k$ as the sums of all products of $k$ distinct variables $x_i$, with negative sign when $k$ is odd, that is
\begin{equation}
  \label{eq:elemsym}
  \tilde{x}_k = (-1)^k \sum_{i_1<i_2<\ldots<i_k} x_{i_1} x_{i_2} \cdots x_{i_k}
\end{equation}
with $\tilde{x}_0 = 1$. The alternating sign $(-1)^k$ in equation~(\ref{eq:elemsym}), which does not appear in the usual definition of the elementary symmetric polynomials, is introduced here as it leads to simpler expressions. While elementary symmetric polynomials are usually denoted using a letter (e.g.\ $\Pi_k$, $s_k$ or $e_k$) which is different from the variable name, we use here the superimposed tilde over the variable name since we will be dealing with elementary symmetric polynomials of different sets of variables.

The \emph{fundamental theorem of symmetric polynomials} states that any symmetric polynomial in the variables $x_i$ can be expressed as a polynomial in the elementary symmetric polynomials $\tilde{x}_k$.

Equation~(\ref{eq:elemsym}) allows computing the elementary symmetric polynomials $\tilde{x}_k$ in terms of the $n$ variables $x_i$. Conversely, the values $x_i$ can be calculated from $\tilde{x}_k$ as the solutions for $x$ of the polynomial equation
\begin{equation}
  \sum_{j=0}^{n} \tilde{x}_{n-j} x^j = 0
\end{equation}

\section{Formulating the system of equations}
\label{sec:formulating}

\subsection{Moment equations}
\label{sec:momeqs}

Our objective is to derive a cubature formula (or \emph{rule}) for the approximate evaluation of the integral of a function $f$ over the area $A$ of a triangle
\begin{equation}
  \label{eq:gencub}
  \bar{I} = \sum_{i=1}^{n_p} \bar{w}_i f^{(i)}
          \approx  \frac{1}{A} \int_A f \,\mathrm{d} A
\end{equation}
where $f^{(i)}$ is the value of $f$ at point $i$, $\bar{w}_i$ is the corresponding weight and $n_p$ is the number of points used in the cubature. We only consider rules of (polynomial) degree $d$, that is rules where equation~(\ref{eq:gencub}) is exact for all polynomials of degree less or equal to $d$, while it is not exact for at least one polynomial of degree $d+1$.

A polynomial of degree $d$ on the triangle can be written as a linear combination of terms $L_1^i L_2^j L_3^{d-i-j}$, where $L_1$, $L_2$ and $L_3$ are the areal coordinates. The cubature rule can therefore be determined by requiring that equation~(\ref{eq:gencub}) is exact for each of these terms. The resulting equations are known as the \emph{moment equations}. The number $\bar{n}_e$ of different terms $L_1^i L_2^j L_3^{d-i-j}$, which is the number of equations to be solved, is
\begin{equation}
  \label{eq:ne_bar}
  \bar{n}_e = (d+1)(d+2)/2
\end{equation}

We only consider \emph{fully symmetric} rules where, if a point with areal coordinates ($L_1$,$L_2$,$L_3$) is used in the cubature, then all points resulting from the permutation of the areal coordinates are also used, with the same weight.
Integration points in a fully symmetric rule can thus belong to one of three different types of point sets, or \emph{orbits}, depending on the number of areal coordinates which are equal. If all areal coordinates are equal, we get a single ``type-0'' orbit, with one point (the centroid). If only two areal coordinates are equal, then we get ``type-1'' orbits with three points which lie on the medians of the triangle. Finally, if all three coordinates are different we get ``type-2'' orbits with six points.
A rule that uses $n_0$ type-0 orbits, $n_1$ type-1 orbits and $n_2$ type-2 orbits is called a rule of type $[n_0,n_1,n_2]$. The number of points for such a rule is
\begin{equation}
  \label{eq:np}
  n_p=n_0 + 3 n_1 + 6 n_2
\end{equation}

Due to the full symmetry employed, when integrating any of the quantities $L_1^i L_2^j L_3^{d-i-j}$ the sum in equation~(\ref{eq:gencub}) will only contain terms of the form
\begin{equation}
   L_1^i L_2^j L_3^{d-i-j} +
   L_1^i L_3^j L_2^{d-i-j} +
   L_2^i L_1^j L_3^{d-i-j} +
   L_2^i L_3^j L_1^{d-i-j} +
   L_3^i L_1^j L_2^{d-i-j} +
   L_3^i L_2^j L_1^{d-i-j}
\end{equation}
These terms are symmetric polynomials, and can therefore be written in terms of the elementary polynomials $\tilde{L}_1 = -(L_1+L_2+L_3)$, $\tilde{L}_2 = L_1 L_2 + L_2 L_3 + L_3 L_1$ and $\tilde{L}_3 = -L_1 L_2 L_3$. It is easily seen that only terms of the form $\tilde{L}_1^k \tilde{L}_2^l \tilde{L}_3^m$ with $k+2l+3m=d$ will be used. Indeed, since $\tilde{L}_1=-1$, only terms of the form $\tilde{L}_2^l \tilde{L}_3^m$ with $2l+3m \leq d$ are actually needed.

The cubature rule of order $d$ can therefore be obtained by requiring that equation~(\ref{eq:gencub}) is exact when the function $f$ is any of the terms $\tilde{L}_2^l \tilde{L}_3^m$ with $2 l + 3 m \leq d$. The number of non-negative solutions of $2 l + 3 m \leq d$ for $l$ and $m$, and therefore the number of equations that must be solved, is given by~\citep{oeisA001399}
\begin{equation}
  \label{eq:numeqs}
  n_e = 1+ \left\lfloor \frac{d^2+6 d}{12} \right\rfloor
\end{equation}
with $\lfloor x \rfloor$ denoting the largest integer that is less or equal to $x$. This is a significant reduction in the number of equations, approximately by a factor of 6 for large values of $d$, compared to the value $\bar{n}_e$ given in equation~(\ref{eq:ne_bar}) for the general case.

While areal coordinates allow for simple formulations of expressions on a generic triangle, they have the disadvantage of introducing three coordinates, instead of the two independent coordinates needed. For this reason, moment equations have generally been obtained using Cartesian or polar coordinates and referring to a specific triangle (exploiting the fact that all triangles are affine). In the fully symmetric case, however, we see that using areal coordinates we easily end up with only two ``coordinates'', the symmetric polynomials $\tilde{L}_2$ and $\tilde{L}_3$.

As will be seen shortly, the moment equations can be further simplified by using, instead of $\tilde{L}_2$ and $\tilde{L}_3$, the quantities
\begin{equation}
  \label{eq:pqdef}
  p = 1-3 \tilde{L}_2 \qquad \text{and}
  \qquad
  q = 1 - \frac{27}{2} \tilde{L}_3 - \frac{9}{2} \tilde{L}_2
\end{equation}
The cubature rule of order $d$ can therefore be obtained by requiring that equation~(\ref{eq:gencub}) is exact when the function $f$ is any of the terms $p^i q^j$ with $2 i + 3 j \leq d$ and $i,j \geq 0$.
The moment equations for a fully symmetric rule of degree~$d$ and type $[n_0,n_1,n_2]$ can thus be written as
\begin{equation}
  \label{eq:eqs6gen}
  \sum_{k=1}^{n_0}   \bar{w}_{0,k} p_{0,k}^i q_{0,k}^j +
  \sum_{k=1}^{n_1} 3 \bar{w}_{1,k} p_{1,k}^i q_{1,k}^j +
  \sum_{k=1}^{n_2} 6 \bar{w}_{2,k} p_{2,k}^i q_{2,k}^j
   = I_{i,j} \quad\text{with}\quad 2i+3j \leq d
\end{equation}
The right hand sides are the integrals
\begin{equation}
  I_{i,j} =  \frac{1}{A} \int_A p^i q^j  \,\mathrm{d} A
\end{equation}
which can be easily computed analytically to give
\begin{equation}
  I_{0,0} = 1,\: I_{1,0} = 1/4,\: I_{0,1} = 1/10,\: I_{2,0} = 1/10,\: I_{1,1} = 2/35,\: I_{3,0} = 29/560,\: I_{0,2} = 7/160,\: I_{2,1} = 1/28, \ldots
\end{equation}

The main advantage of using the quantities $p$ and $q$ is that for type-1 orbits we can introduce a new variable $u$ so that $p=u^2$ and $q=u^3$ and therefore $p^i q^j = u^{2 i + 3 j}$, while for the type-0 orbit $p=q=0$. Setting $w_0 = \bar{w}_{0,1}$, $v_k = 3 \bar{w}_{1,k}$ and $w_k = 6 \bar{w}_{2,k}$, after some computations, the moment equations are finally written as
\begin{subequations}
\label{eq:momeqs}
\begin{gather}
  \label{eq:momeq1}
  w_0 + \sum_{k=1}^{n_1} v_k + \sum_{k=1}^{n_2} w_k = I_{0,0}\\
  \label{eq:momeq3}
  \sum_{k=1}^{n_1} v_k u_k^{2i+3j} + \sum_{k=1}^{n_2} w_k p_k^i q_k^j = I_{i,j}
  \quad\text{with}\quad  2i+3j \leq d, j \leq 1\\
  \label{eq:momeq6}
  \sum_{k=1}^{n_2} w_k (p_k^3-q_k^2) p_k^i q_k^j = I_{i+3,j} - I_{i,j+2}
  \quad\text{with}\quad 2i+3j \leq d-6
\end{gather}
\end{subequations}
where in equation~(\ref{eq:momeq1}) we set $w_0=0$ if $n_0=0$.

For both $d=0$ and $d=1$ the only moment equation is~(\ref{eq:momeq1}). This means that any (fully symmetric) rule exact for $d=0$ will also be exact for $d=1$, thus there are no rules of degree 0. For this reason in the following we always assume that $d \geq 1$.

\subsection{Consistency conditions}
\label{sec:conscond}

To set up the moment of equations for a rule of degree $d$, it is first necessary to determine the type of the rule, i.e. the number of orbits of each type.

The moment equations~(\ref{eq:momeqs}) form a system of $n_e$ equations in $n_v$ variables, where $n_e$ is given in equation~(\ref{eq:numeqs}) while $n_v = n_0 + 2 n_1 + 3 n_2$. Similarly, the subsystem~(\ref{eq:momeq6}) has $n_e-d$ equations and $3 n_2$ variables.

We assume that both the system~(\ref{eq:momeqs}) and its subsystem~(\ref{eq:momeq6}) are inconsistent (i.e.\ have no solutions) if and only if they are overdetermined (i.e. have more equations than variables). This assumption, together with the fact that there may be at most one type~0 orbit, yields the following \emph{consistency conditions}
\begin{subequations}
\begin{gather}
  3 n_2 \geq n_e - d \\
  3 n_2 + 2 n_1 + n_0 \geq n_e \\
  n_0 \leq 1
\end{gather}
\end{subequations}
which must be satisfied to obtain a solution of the moment equations, and thus they restrict the choice of the rule type. For a given degree $d$, a minimal-point rule is sought, that is a rule that satisfies the consistency conditions with the lowest total number of points, as given by equation~(\ref{eq:np}). This yields
\begin{equation}
  n_2 = \big\lfloor (n_e - d + 2)/3 \big\rfloor
  ,\qquad
  n_1 = \big\lfloor (n_e - 3 n_2)/2 \big\rfloor
  ,\qquad
  n_0 = n_e - 3 n_2 - 2 n_1
\end{equation}

It is conceivable that a rule that violates the consistency conditions may lead to a system of moment equations that, although overdetermined, has solutions. These so-called \emph{fortuitous} rules have great theoretical interest, as well as practical interest in the case where they have fewer integration points compared to the minimal-point rules described above. No fortuitous rules are encountered in the present paper, however, nor in the available literature on cubature rules on the triangle.

The system of moment equations~(\ref{eq:momeqs}) can be inconsistent, zero-dimensional or positive-dimensional (with zero solutions, a finite number of solutions or infinite solutions respectively). We use here the same terms to identify the corresponding rule types and individual rules, thus we have inconsistent rule types, which yield no rules, zero-dimensional rule types, which yield a finite number of zero-dimensional rules, and  positive-dimensional rule types which yield an infinite number of positive dimensional rules. In the case of positive-dimensional rule types, the analytical solution can be expressed using a number of free parameters.

\subsection{Advantages of the suggested form of the moment equations}

The development of the method given in Sections~\ref{sec:momeqs} and~\ref{sec:conscond} to formulate the moment equations using symmetric polynomials follows in some main points the classic one presented by \citet{LynessJespersen1975}. It has, however, the obvious benefit of providing polynomial moment equations, while~\citep{LynessJespersen1975} also uses cosines. In this, the present method is similar to the one presented by \citet{WandzuraXiao2003}.

All three methods  are equivalent, in that they yield the same rules. Indeed, it is relatively easy to pass from one method to the other: setting $p_i = r_i^2$, $q_i = r_i^3 \cos 3 \alpha_i$ and $u_i = r_i$ in equations~(\ref{eq:momeqs}) yields the moment equations in~\citep{LynessJespersen1975}, while it is easily seen that, for the triangle used in~\citep{WandzuraXiao2003}, $p$ and $q$ are equal to the invariants $x^2+y^2$ and $x^3-3 x y^2$.

The present method is arguably simpler and more intuitive in its formulation, while it provides simpler formulas. Additionally, this method is elegantly formulated without reference to a specific triangle. From a practical point of view, however, the main advantage is that the resulting polynomial equations are of significantly lower degree than those provided by the other methods, for example the maximum degree of equations~(\ref{eq:momeq6}) is $\lfloor d/2 \rfloor + 1$ instead of $d+1$. This is especially important when solving the equations analytically.

\section{Analytical solution of the moment equations}
\label{sec:analytsol}

\subsection{The usefulness of analytical solutions}

Except for some trivial low-degree rules, the moment equations are generally solved numerically, e.g.\ using a multivariate Newton-Raphson solver. The cubature rule is then given as a table of integration point coordinates and weights, expressed as floating point approximations of a given precision. This numerical approximation of the cubature rule is the one actually required when using the rule in applications.

Numerical methods have the advantage of being able to provide cubature rules of high degree~\citep[see e.g.][]{XiaoGimbutas2010}. Convergence of the method to a solution is not guaranteed, however, as it most often depends on the selection of an appropriate ``initial guess'' required by the solver. This means that inability to obtain a solution does not prove that the solution does not exist. Additionally, when a solution is obtained numerically, no information is obtained regarding the existence of other solutions.
For this reason, in this paper we investigate the analytical solution of the moment equations, in order to obtain a definitive answer regarding the different cubature rules for a given degree and type.

There exist algorithms for solving analytically arbitrary systems of polynomial equations, for example using Gr\"obner bases (see~\citep{Lazard2009} for an informal overview of the state of the art). Unfortunately, their requirements in both computer memory and computation time are such that in practice they fail to provide a solution even for rules of relatively low degree. To obtain solutions for higher degrees, it is therefore necessary to exploit as much as possible the structure of the moment equations.

An interesting alternative to the analytical solution of the moment equations is to use homotopy continuation methods to compute \emph{all} solutions of the system numerically~\citep{Verschelde1999}. This is however clearly beyond the scope of the present paper.

\subsection{Solution strategy}

The subsystems~(\ref{eq:momeq1}), (\ref{eq:momeq3}) and~(\ref{eq:momeq6}) have respectively $1$, $d-1$ and $n_e-d$ equations. The weight $w_0$ (if it is non-zero) appears only in equation~(\ref{eq:momeq1}) while the variables $v_k$ and $u_k$ appear only in equations~(\ref{eq:momeq1}) and (\ref{eq:momeq3}).

Consider first the case of a rule with a type-0 orbit ($n_0=1$). Equation~(\ref{eq:momeq1}) is then just used to determine $w_0$ when all other weights have been calculated. The weights $v_k$ of type-1 orbits can be eliminated from equations~(\ref{eq:momeq3}), as described in~\citep[pp.~771--773]{RabinowitzRichter1969} for cubature rules on other regions, to obtain the (linear in the symmetric polynomials $\tilde{u}_k$) system of equations
\begin{equation}
  \label{eq:momeq3sym}
  \sum_{k=0}^{n_1} J_{i-k} \tilde{u}_k = 0, \qquad i=n_1+2, \ldots, d
\end{equation}
where
\begin{equation}
  \label{eq:Jdef}
  J_i = \begin{cases}
    \displaystyle{ I_{j,0} - \sum_{k=1}^{n_2} w_k p_k^j     } & \text{if $i=2j$}   \\[2em]
    \displaystyle{ I_{j,1} - \sum_{k=1}^{n_2} w_k p_k^j q_k } & \text{if $i=2j+3$}
  \end{cases}
\end{equation}

The system~(\ref{eq:momeq3sym}) has $n_1$ unknowns $\tilde{u}_k$ (since $\tilde{u}_0 = 1$) and $d-n_1-1$ equations. If $n_1=(d-1)/2$ then equations~(\ref{eq:momeq6}) are sufficient to evaluate the variables $w_k$, $p_k$ and $q_k$ of type-2 orbits, and then equations~(\ref{eq:momeq3sym}), (\ref{eq:momeq3}) and~(\ref{eq:momeq1}) yield in turn the values of $\tilde{u}_k$, $v_k$ and $w_0$. The same happens if $n_1>(d-1)/2$, but in this case the system is positive-dimensional and some of the $\tilde{u}_k$ remain as free parameters in the solution. Finally, if $n_1<(d-1)/2$ then obtaining a solution is more difficult, since to evaluate $w_k$, $p_k$ and $q_k$ we need not only equations~(\ref{eq:momeq6}) but also the equations that remain after eliminating $\tilde{u}_k$ from~(\ref{eq:momeq3sym}).

When the type-0 orbit is not used ($n_0=0$), it is generally easier to introduce an additional equation
\begin{equation}
  \label{eq:momeqextra}
  \sum_{k=1}^{n_1} v_k u_k = J_1
\end{equation}
where $J_1$ is an unknown quantity, which is not defined by~(\ref{eq:Jdef}). Eliminating the weights $v_k$ from equations~(\ref{eq:momeq1}), (\ref{eq:momeq3}) and~(\ref{eq:momeqextra}) leads to a system of equations like~(\ref{eq:momeq3sym}), only that the index $i$ is now in the range $i=n_1, \ldots, d$ and $J_1$ is an additional unknown that must be eliminated.

In all cases, equations~(\ref{eq:momeq6}) must be solved, possibly together with the equations that remain after eliminating $\tilde{u}_k$ from~(\ref{eq:momeq3sym}). Unfortunately, no easy way has been found to simplify these equations as we did to derive the system~(\ref{eq:momeq3sym}). The use of symmetric polynomials can, however, again lead to somehow simpler expressions.

\subsection{Permutation invariance of the orbits}
\label{sec:orbinv}

In Section~\ref{sec:formulating} we exploited the fact that the cubature rule is invariant with respect to a permutation of the integration points within a given orbit, and expressed this invariance using symmetric polynomials.

Another obvious property of the cubature rules, which however has received much less attention in the literature, is their invariance with respect to permutation of orbits of the same type. This is reflected in the fact that the moment equations~(\ref{eq:momeqs}) are polynomials which are ``symmetric'' (i.e.\ invariant with respect to permutation) in the pairs ($u_k$,$v_k$) and in the triplets ($p_k$, $q_k$, $w_k$). This can be seen from the system~(\ref{eq:momeq3sym}) where, having eliminated the $v_k$, the resulting polynomials are symmetric in the $u_k$ and have thus been expressed in terms of the elementary symmetric polynomials $\tilde{u}_k$. In a similar way, eliminating $q_k$ and $w_k$ allows us to express the moment equations in terms of the symmetric polynomials $\tilde{p}_k$.

The system that results by eliminating the $v_k$, $q_k$ and $w_k$ from the moment equations and expressing the results in terms of the $\tilde{u}_k$ and $\tilde{p}_k$ is in most cases more complicated than the moment equations. It has however fewer variables, and it leads to a much simpler expression for the solution, when such a solution is actually found.

Indeed, one important advantage of expressing the moment equations in terms of symmetric polynomials is that the number of solutions of the system is equal to the number of different cubature rules that can be obtained. Consider for example the degree-4 $[0,2,0]$ rule, for which~\citet{LynessJespersen1975} mention that, in the present notation, $u_1$ and $u_2$ are the roots of $15 x^4 + 20 x^3 - 30 x^2 + 4$. This does not mean, however, than any combination of the roots is a valid solution for $u_1$ and $u_2$, indeed only two pairs of solutions give a cubature rule. In terms of symmetric polynomials, on the other hand, the solution is obtained by solving the equations $3\tilde{u}_1^2-4\tilde{u}_1-2 = 0$ and $5\tilde{u}_2+2\tilde{u}_1+2 = 0$, where it is seen that two different rules are obtained, one for each solution of the system.

It is worth considering that even when solving the moment equations numerically, considering the invariance with respect to permutation of orbits of the same type can have a significant effect on the solution method. As an example, there is only one degree-15 [1,7,4] rule. The system~(\ref{eq:momeq6}) however has $4!=24$ solutions, while if we were to solve all equations~(\ref{eq:momeqs}) together we would have $7!4!=120960$ solutions. It is thus conceivable that an iterative numerical solution algorithm may fail to converge by being ``attracted'' in turn by different solutions.

\subsection{Solution quality}

Once a cubature rule is determined by solving the moment equations, the sign of the weights and the position of the integration points is examined, to determine the \emph{quality} of the solution. The quality is described using a two-letter label: the first letter is P if all weights are positive and N if at least one weight is negative, while the second letter is I if all points are inside the triangle, O if there is at least one point outside the triangle, and B if no points are outside the triangle but at least one is on the boundary of the triangle. The following qualities are therefore encountered: PI, NI, PB, NB, PO, NO.

In all the above cases, the coordinates and weights of the integration points are considered to be real. Though it is well-known that complex solutions may exist, these are not taken into account, since a cubature rule with complex-valued coordinates of the integration points would be of little, if any, use. Moreover, the moment equations are usually solved using numerical methods that only return real solutions, as these methods perform significantly better than methods that could return complex solutions.

On the other hand, when obtaining the solutions analytically it costs nothing to also consider complex solutions. For this reason, we expand the above definition of the quality of cubature rules by setting the first letter of the label to C if at least one weight is complex-valued and by setting the second letter of the label to C if at least one integration point has complex coordinates. Interestingly, while it is not possible to have complex weights without complex coordinates, it is possible to have real weights with complex coordinates. The following three additional qualities are therefore obtained: CC, PC, NC.

Including complex solutions allows us to make the distinction between moment equations that have no solution and those that have solutions, even though they may all be complex. Considering as an example a degree-15 rule, there are no solutions for type $[0,7,4]$ (which does not satisfy the consistency conditions), while there is a single complex (NC) solution for type $[1,7,4]$ (which satisfies the consistency conditions). It is generally  expected that all types satisfying the consistency conditions will yield at least one solution, but with complex solutions appearing with increasing frequency as the degree of the rule increases.

Although we compute all solutions, independently of their quality, in most applications we need rules of PI (or at most NI) quality. For this reason, if a minimal-point rule does not yield any PI rules, we investigate rules with increasingly more points until a rule is found that has a PI solution. When considering rules with additional points, it is possible to have rules with the same degree and number of points, but different type and different number of free parameters appearing in the solution.

Consider for example the degree-7 rules. The minimal-point rule $[1,2,1]$ has 13 points and the best quality achievable with it is NI. Increasing the number of points, we get either a $[0,3,1]$ or a $[0,1,2]$ rule, both with 15 points, where the first has one free parameter while the second has none. In this case, where both types can yield PI rules, we prefer the zero-dimensional one as it has more type-2 orbits, so less integration points are restricted to be located on the medians.

In general, among rules with the same number of points and the same quality, we would prefer those with more type-2 orbits and thus less free parameters. The presence of free parameters in the solution of the moment equations, on the other hand, allows for much greater flexibility in obtaining a rule of PI quality. Moreover, the use of more type-1 orbits leads to simpler moment equations, which are easier to solve analytically.

Note that the numerical, iterative solution of the moment equations for positive-dimensional rules \citep[see e.g.][]{WandzuraXiao2003} yields only one of the infinite solutions. Though it is possible to consider numerically the variation of the solution with the variation of a parameter \citep[see e.g.][]{BerntsenEspelid1990}, analytical solutions are much more powerful in studying parametrically positive-dimensional cubature rules and their quality. The study and presentation of such rules, however, requires a much more extensive discussion which goes well beyond the scope of the present paper. For this reason, in Section~\ref{sec:results} we only present results for zero-dimensional cubature rules.

\section{Results and discussion}
\label{sec:results}

Using the method described in Sections~\ref{sec:formulating} and~\ref{sec:analytsol} we compute here analytically cubature rules for degree up to~15. As described in Section~\ref{sec:orbinv}, the permutation invariance of the orbits should be exploited to express the moment equations~(\ref{eq:momeqs}) in a form more suitable for analytical solution, for example in terms of the symmetric polynomials $\tilde{u}_k$ and $\tilde{p}_k$. This has been achieved for each degree and rule type in a heuristic way, which involved (for higher degrees) extensive calculations until the initial system was transformed into a new one, solvable (on the available hardware and software) using Gr\"obner bases. The actual calculations performed in each case are obviously too lengthy to be written out here. Indeed, in the non-trivial cases, the analytical solution itself becomes too long, as is already apparent in Appendix~\ref{sec:app_analyt} for the degree-6 rule.

Table~\ref{tab:rulesummary} gives a summary of the properties of all cubature rules thus computed. As already mentioned, we only consider zero-dimensional rules. We calculate for each degree the minimal-point rules and, if none of these are of quality PI, we calculate additional rule types with more points until a rule with PI quality is found (except for $d=15$ where additional rules were not computed). Appendix~\ref{sec:app_analyt} provides analytical expressions for evaluating some of the cubature rules, while Appendix~\ref{sec:app_numer} provides numerical values for new rules of PI or NI quality.

\begin{table}
\centering
\caption{Summary of the properties of all computed rules}\label{tab:rulesummary}
\begin{tabular}{rcrrrrrrrrrr}
\toprule
degree & type & points & solutions & PI & NI & PB & PO & NO & PC & NC & CC \\
\midrule
  1 & $[1,0,0]$ &  1 &  1   &  1 & -- & -- & -- & -- & -- & -- & -- \\
  2 & $[0,1,0]$ &  3 &  2   &  1 & -- &  1 & -- & -- & -- & -- & -- \\
  3 & $[1,1,0]$ &  4 &  1   & -- &  1 & -- & -- & -- & -- & -- & -- \\
    & $[0,0,1]$ &  6 &  1   &  1 & -- & -- & -- & -- & -- & -- & -- \\
  4 & $[0,2,0]$ &  6 &  2   &  1 & -- & -- &  1 & -- & -- & -- & -- \\
  5 & $[1,2,0]$ &  7 &  1   &  1 & -- & -- & -- & -- & -- & -- & -- \\
  6 & $[0,2,1]$ & 12 &  6   &  2 & -- & -- &  2 & -- & -- & -- &  2 \\
  7 & $[1,2,1]$ & 13 &  4   & -- &  1 & -- &  1 & -- & -- & -- &  2 \\
    & $[0,1,2]$ & 15 &  4   &  2 & -- & -- & -- & -- & -- & -- &  2 \\
  8 & $[1,3,1]$ & 16 &  2   &  1 &  1 & -- & -- & -- & -- & -- & -- \\
  9 & $[1,4,1]$ & 19 &  1   &  1 & -- & -- & -- & -- & -- & -- & -- \\
 10 & $[0,4,2]$ & 24 & 14   & -- & -- & -- &  4 &  1 & -- & -- &  9 \\
    & $[1,2,3]$ & 25 & 15   &  4 & -- & -- & -- &  2 & -- &  3 &  6 \\
 11 & $[0,5,2]$ & 27 &  6   & -- & -- & -- &  1 & -- & -- &  2 &  3 \\
    & $[1,3,3]$ & 28 & 23   & -- &  2 & -- &  5 &  3 &  2 &  4 &  7 \\
    & $[0,2,4]$ & 30 & 34   &  4 & -- & -- &  1 &  1 &  4 &  2 & 22 \\
 12 & $[0,5,3]$ & 33 & 24   &  2 &  1 & -- & -- & -- & -- & -- & 21 \\
 13 & $[0,6,3]$ & 36 &  8   & -- & -- & -- & -- &  1 & -- &  1 &  6 \\
    & $[1,4,4]$ & 37 & 54   &  2 &  3 & -- &  4 &  5 &  2 &  8 & 30 \\
 14 & $[0,6,4]$ & 42 & 38   &  1 & -- & -- &  3 &  3 & -- &  3 & 28 \\
 15 & $[1,7,4]$ & 46 &  1   & -- & -- & -- & -- & -- & -- &  1 & -- \\
\bottomrule
\end{tabular}
\end{table}

The only case where three rule types must be computed to obtain PI quality is $d=11$. This is therefore the only case (for $d<15$) where a positive-dimensional rule of PI quality (type $[1,5,2]$ with 28 points) has less points than the best possible zero-dimensional rule of the same quality (type $[0,2,4]$ with 30 points).
\footnote{Type $[1,5,2]$ fully symmetric PI rules exist, and are easy to obtain using the method presented in this paper, yet no such rule was encountered in the literature. \citet{LynessJespersen1975} present a $[1,5,2]$ PB rule, \citet{WandzuraXiao2003} compute what is most probably a $[1,5,2]$ PI rule but do not present it, while the 28-point PB rule given by \Citet{Taylor2007} and the 27-point PI rule given by \citet{Taylor2008} are asymmetric rules.}
Where NI rules are acceptable, the two $[1,3,3]$ rules can be used, since they have the same number of points as the $[1,5,2]$ rules and have less points on the medians. The second $[1,3,3]$ NI rule given in Appendix~\ref{sec:app_numer} is then to be preferred as it has a small negative weight for a single point, while the first one has a large negative weight for three points.

The results summarised in Table~\ref{tab:rulesummary} confirm the general expectation that as the rule degree increases the number of solutions will increase, though with most solutions being complex ones. This is not always the case, however, as evidenced by the existence of a single $[1,7,4]$ rule for $d=15$. It is thus clear that it is not possible to detect in these results a specific pattern in the number of solutions, the number of real solutions or the number of PI (or NI) solutions.

An interesting side effect of computing analytically the cubature rules is that we can prove the non-existence of specific fortuitous rules. Consider for example the case of degree-10 rules where we compute the 24-point $[0,4,2]$ rule. If there were a fortuitous rule with two type-2 orbits and less than 24 points, then the $[0,4,2]$ rule should be positive-dimensional in order to depend on some parameters which, for specific values, would yield the fortuitous rule. Computing analytically the $[0,4,2]$ rule, however, shows that it is zero-dimensional, as expected. Similarly, since the $[1,2,3]$ rule is zero-dimensional, there exist no fortuitous rules with three type-2 orbits and less than 25 points. Since it is easily shown that for degree 10 no rules exist with one or zero type-2 orbits and also that no $[0,0,4]$ rules exist (which would have 24 points) we see that there are no fortuitous degree-10 rules with 24 points or less. Similar tests can be performed for all other rule degrees considered here.

A list of numerical values for all computed rules, independently of their quality, is provided together with this paper as supplemental material. This list includes the zero-dimensional rules found in~\citep{Stroud1971, Cowper1973, LynessJespersen1975, LaursenGellert1978, Dunavant1985}. An interesting property of some rules of bad quality (i.e. neither PI nor NI) is that the orbits that have points outside the triangle or with complex coordinates have a much smaller weight (in absolute value). This is the case for example for $d=11$ and the fourth $[1,3,3]$ NC rule, or for $d=15$ and the $[1,7,4]$ NC rule. These rules, together with a node elimination algorithm \citep{XiaoGimbutas2010}, could possibly be used to derive   cubature rules that are not fully symmetric with fewer points than the fully symmetric ones.

\section{Conclusions}

In this paper we have used symmetric polynomials to express the double invariance inherent in fully symmetric cubature rules in the triangle (invariance with respect to permutation of points within an orbit and with respect to permutation of orbits of the same type). This has allowed us to formulate the moment equations in such a way that analytical solutions have been derived for zero-dimensional rules of degree up to 15.

A few new rules of good quality have been thus derived and are given in Appendix~\ref{sec:app_numer}. Additionally, the analytical solutions ensure that all possible rules of a given type and degree were computed, independently of their quality. This allows us, for example, to prove that indeed no rules of good quality (PI or even NI) exist for some cases where no such rules were encountered in the literature.

Though only zero-dimensional rules have been computed here, the proposed analytical approach is also well-suited for the thorough study of positive-dimensional rules. In this case, however, an additional difficulty lies in finding intuitive and useful ways to present the (infinite) solutions and their properties.

In all cases, combining a better understanding of the structure of the moment equations together with better-performing algorithms, software implementations and hardware platforms, should allow determining rules of increasingly high degree.

\section*{Acknowledgements}

The research leading to these results has received funding from the
European Research Council under the European Community's Seventh
Framework Programme (FP7/2007--2013) / ERC grant agreement n\textsuperscript{o} 228051 [MEDIGRA].

\appendix
\section{Analytical expressions for the cubature rules}
\label{sec:app_analyt}

This appendix lists analytical expressions for some of the rule types considered in this paper. For each rule type we list the degree $d$, the rule type and a list of expressions.  Using these expressions, it is easy to obtain the coordinates and weights of the integration points for all rules of the given type.

\begin{itemize}
\item $d=1$\quad $[1,0,0]$:\quad $w_0=1$

\item $d=2$\quad $[0,1,0]$:\quad $u_1^2 = 1/4$, $v_1=1$

\item $d=3$\quad $[1,1,0]$:\quad $u_1 = 2/5$, $v_1=25/16$, $w_0=-9/16$

\item $d=3$\quad $[0,0,1]$:\quad $p_1 = 1/4, q_1 = 1/10, w_1 = 1$

\item $d=4$\quad $[0,2,0]$:\\
      $3 \tilde{u}_1^2 - 4 \tilde{u}_1 - 2 = 0$, $\tilde{u}_2=-2/5 (\tilde{u}_1+1)$,
      $v_i = \big((81/248)\tilde{u}_1-6/31\big) u_i + (15/124) \tilde{u}_1 + 151/248$

\item $d=5$\quad $[1,2,0]$: $\tilde{u}_1 = -2/7$, $\tilde{u}_2= -2/7$, $v_i = -(7/400) u_i + 39/100$, $w_0=9/40$

\item $d=6$\quad $[0,2,1]$:\\[0.1em]
$p_1^6 - \frac{2943}{896} p_1^5 + \frac{12577377}{3211264} p_1^4 - \frac{6335029}{2809856} p_1^3 + \frac{211997025}{314703872} p_1^2 - \frac{7914723}{78675968} p_1 + \frac{14953009}{2517630976} = 0$,\\[0.2em]
$q_1= \frac{6773849}{1180960} - \frac{17597477}{258335} p_1 + \frac{618894079}{2066680} p_1^2 - \frac{2009158}{3355} p_1^3 + \frac{20120576}{36905} p_1^4 - \frac{6422528}{36905} p_1^5$,\\[0.2em]
$w_1=
 \frac{88271353265388941906672}{1552328339949698669325}
-\frac{2093018886005051378041487}{3104656679899397338650} p_1
+\frac{4677969412268735483683874}{1552328339949698669325} p_1^2$\\[0.2em]\hspace*{2em}$
-\frac{874483029603676756153618}{141120758177245333575}   p_1^3
+\frac{8938712246012353125723136}{1552328339949698669325} p_1^4
-\frac{2885760563751222732259328}{1552328339949698669325} p_1^5$,\\[0.3em]
$\tilde{u}_1 =
-\frac{5647577278829}{5843759130}
+\frac{219725386019839}{17838843660} p_1
-\frac{4934508553334726}{84734507385} p_1^2$\\[0.2em]\hspace*{2em}$
+\frac{965776421126167}{7703137035} p_1^3
-\frac{538505362157056}{4459710915} p_1^4
+\frac{3379769853673472}{84734507385} p_1^5$,\\[0.2em]
$\tilde{u}_2 =
 \frac{15104616525664}{20453156955}
-\frac{581971572152849}{62435952810} p_1
+\frac{3683852401439816}{84734507385} p_1^2$\\[0.2em]\hspace*{2em}$
-\frac{709365733908202}{7703137035} p_1^3
+\frac{389416992756736}{4459710915} p_1^4
-\frac{2417750126231552}{84734507385} p_1^{5}$,\\[0.2em]
$v_i = \big(
 \frac{409434268039529549940720811615256576}{158205456880303475487279097255725} p_1^5
-\frac{410497105230467053184700451899528704}{52735152293434491829093032418575} p_1^4
$\\[0.2em]\hspace*{2em}$
+\frac{114666485932207310484775951381500251}{14382314261845770498843554295975} p_1^3
-\frac{1141034063408380347314793772529581321}{316410913760606950974558194511450} p_1^2
$\\[0.2em]\hspace*{2em}$
+\frac{153804297816157841896911608744616331}{210940609173737967316372129674300} p_1
-\frac{68148358857994347902974327222077377}{1265643655042427803898232778045800}
\big) u_i
$\\[0.2em]\hspace*{2em}$
-\frac{660366862842903249663697383981056}{606151175786603354357391177225} p_1^5
+\frac{17884598780372115712297691281128448}{5455360582079430189216520595025} p_1^4
$\\[0.2em]\hspace*{2em}$
-\frac{1660804293851424897302836415981834}{495941871098130017201501872275} p_1^3
+\frac{2716663212121457459566848039780239}{1818453527359810063072173531675} p_1^2
$\\[0.2em]\hspace*{2em}$
-\frac{3169975113118311937146770957038031}{10910721164158860378433041190050} p_1
+\frac{220460485921384140338311776720617}{10910721164158860378433041190050}$

\item $d=7$\quad $[1,2,1]$:\\[0.1em]
$p_1^4 - \frac{23}{12} p_1^3 + \frac{655}{448} p_1^2 - \frac{85}{196} p_1 + \frac{1619}{37632} = 0$, \\[0.2em]
$q_1=\frac{73}{160} - \frac{63}{20} p_1 + \frac{273}{40} p_1^2-\frac{21}{5} p_1^3$,
$w_1= \frac{5559373039}{1374543450} - \frac{4035503891}{196363350} p_1 +  \frac{3029805464}{98181675} p_1^2 - \frac{577446688}{32727225} p_1^3$,\\[0.2em]
$\tilde{u}_1 = \frac{204779}{4630} - \frac{616196}{2315} p_1 + \frac{978558}{2315} p_1^2 - \frac{585648}{2315} p_1^3$,
$\tilde{u}_2 = -\frac{86623}{2315} + \frac{511579}{2315} p_1 - \frac{811132}{2315} p_1^2 + \frac{484512}{2315} p_1^3$,\\[0.2em]
$v_i=
\Big(
-\frac{637366793665532978264}{9052562883613960471} p_1^3
+\frac{3156091037460298906045}{27157688650841881413} p_1^2
-\frac{9803429487627684799252}{135788443254209407065} p_1
+\frac{14666951220214040085227}{1267358803705954465940}
\Big) u_i
$\\[0.2em]\hspace*{2em}$
+\frac{7823399076093706515424}{135788443254209407065} p_1^3
-\frac{39124895515170463542614}{407365329762628221195} p_1^2
+\frac{98256377831808794616331}{1629461319050512884780} p_1
-\frac{109198370776069008639239}{11406229233353590193460}
$,\\[0.2em]
$w_0 = -\frac{3660769728}{100486445} p_1^3 + \frac{4347049032}{703405115} + \frac{6057843876}{100486445} p_1^2 - \frac{752902776}{20097289} p_1$

\item $d=7$\quad $[0,1,2]$:\\[0.1em]
$u_1^4 - (4/9) u_1^3 - (1/3) u_1^2 + (1/36)$,
$v_1=\frac{156673}{8817780} + \frac{2159752}{2204445} u_1 + \frac{5368006}{2204445} u_1^2 - \frac{3133452}{734815} u_1^3$,\\[0.2em]
$\tilde{p}_1=-\frac{1079}{1281} - \frac{310}{1281} u_1 - \frac{128}{427} u_1^2 + \frac{720}{427} u_1^3$,
$\tilde{p}_2 = \frac{1493}{11956} + \frac{1130}{8967} u_1 + \frac{2116}{8967} u_1^2 - \frac{2046}{2989} u_1^3$,\\[0.2em]
$q_i=
\Big(
\frac{489}{427} + \frac{465}{854} u_1 + \frac{288}{427} u_1^2 -\frac{1620}{427} u_1^3 \Big) p_i
-\frac{653}{2989} - \frac{1695}{5978} u_1 - \frac{1587}{2989} u_1^2 + \frac{9207}{5978} u_1^3$,
$w_i=
\Big(
-\frac{21117033567}{4098798070} u_1^3
+\frac{34215330023}{8197596140} u_1^2
+\frac{676519529}{409879807} u_1
-\frac{22884360891}{16395192280}
\Big) p_i
+\frac{7377613908}{2049399035} u_1^3
-\frac{123702006967}{49185576840} u_1^2
-\frac{10640567105}{9837115368} u_1
+\frac{103431908839}{98371153680}
$

\item $d=8$\quad $[1,3,1]$:\\[0.1em]
$p_1=2/5$,
$q_1^2 - \frac{116}{355} q_1 + \frac{443}{17750}$,
$w_1= \frac{1286875}{529326} q_1 - \frac{561275}{4234608}$,
$w_{{0}}={\frac{197671347}{256973920}}-{\frac{32981985}{6424348}} q_1$,\\[0.2em]
$\tilde{u}_1= \frac{181760}{50289} q_1 - \frac{70340}{50289}$,
$\tilde{u}_2=-\frac{124960}{50289} q_1 + \frac{10642}{50289}$,
$\tilde{u}_3=-\frac{50410}{50289} q_1 + \frac{14008}{50289}$,\\[0.2em]
$v_i=\Big(
-\frac{18610928498796911607672845}{2275833709245992090597766} u_i^2
+\frac{19252524428004364259122223}{4551667418491984181195532} u_i
+\frac{12430296145585635931273841}{4551667418491984181195532} \Big) q_1
$\\[0.2em]\hspace*{2em}$
+\frac{11189621308192975569101651}{22758337092459920905977660} u_i^2
-\frac{2452756844382101152643719}{5689584273114980226494415} u_i
+\frac{10275755611647081695669293}{182066696739679367247821280}$

\item $d=9$\quad $[1,4,1]$:\\[0.1em]
$p_1=2/5$,
$q_1=2/11$,
$w_1=\frac{3025}{11648}$,
$w_0=\frac{85293}{878080}$
$\tilde{u}_1=-\frac{212}{407}$,
$\tilde{u}_2=-\frac{1002}{2035}$,
$\tilde{u}_3= \frac{212}{2035}$,
$\tilde{u}_4= \frac{112}{2035}$,\\[0.2em]
$v_i=
 \frac{506023048885425107}{1503746382262924800}
+\frac{11465050245708013}{334165862725094400} u_i
-\frac{10064998401780383}{12531219852191040}  u_i^2
+\frac{52676213406614851}{109363373255485440} u_i^3$

\end{itemize}

\section{Numerical values for new cubature rules of good quality}
\label{sec:app_numer}

This appendix lists all the computed cubature rules of quality PI and NI that are \emph{not} listed in the Encyclopedia of Cubature Formulas~\citep{Cools2003}. For each rule we first list the degree $d$, the number of points $n_p$, the rule type and the rule quality. We then provide a list of the orbits, where the first column is the number of points in the orbit, the second is the weight for each integration point and the last three columns are the areal coordinates defining a point in the orbit.

\begin{center}
\small
\ttfamily
\begin{tabular}{rrrrr}
\multicolumn{5}{c}{\rmfamily $d=7$, $n_p=15$, type $[0,1,2]$, quality PI}\\
\toprule
3 & ~1.253936074493031e-01 & 5.134817203287849e-01 & 2.432591398356075e-01 & 2.432591398356075e-01\\
6 & ~7.630633834054171e-02 & 5.071438430720704e-02 & 3.186441898475371e-01 & 6.306414258452559e-01\\
6 & ~2.766352460147343e-02 & 4.572082984632032e-02 & 8.663663134174900e-02 & 8.676425388119307e-01\\
\bottomrule
\end{tabular}
\bigskip

\begin{tabular}{rrrrr}
\multicolumn{5}{c}{\rmfamily $d=10$, $n_p=25$, type $[1,2,3]$, quality PI}\\
\toprule
1 & ~8.321973698645014e-02 & 3.333333333333333e-01 & 3.333333333333333e-01 & 3.333333333333333e-01\\
3 & ~5.265194946824459e-02 & 6.741737642518105e-01 & 1.629131178740948e-01 & 1.629131178740948e-01\\
3 & ~1.095128834026841e-02 & 9.429929994232243e-01 & 2.850350028838784e-02 & 2.850350028838784e-02\\
6 & ~5.627727971081118e-02 & 1.468115053939304e-01 & 3.366958752782316e-01 & 5.164926193278379e-01\\
6 & ~3.539494779153839e-02 & 2.930760450457947e-02 & 3.633626169945705e-01 & 6.073297785008500e-01\\
6 & ~2.932286409565224e-02 & 3.368569868061029e-02 & 1.533030551695614e-01 & 8.130112461498283e-01\\
\bottomrule
\end{tabular}
\bigskip

\begin{tabular}{rrrrr}
\multicolumn{5}{c}{\rmfamily $d=11$, $n_p=28$, type $[1,3,3]$, quality NI}\\
\toprule
1 & ~1.918874890144834e-01 & 3.333333333333333e-01 & 3.333333333333333e-01 & 3.333333333333333e-01\\
3 & ~4.494673641886435e-02 & 3.618633512713571e-02 & 4.819068324364321e-01 & 4.819068324364321e-01\\
3 & ~4.110012887519764e-02 & 8.110470430369368e-01 & 9.447647848153162e-02 & 9.447647848153162e-02\\
3 & -3.968990834546128e+00 & 5.167437799769402e-01 & 2.416281100115299e-01 & 2.416281100115299e-01\\
6 & ~2.034937507030261e+00 & 2.273118971929394e-01 & 2.539744293918833e-01 & 5.187136734151773e-01\\
6 & ~3.401519357474532e-02 & 3.065503892384041e-02 & 2.512520556604245e-01 & 7.180929054157351e-01\\
6 & ~7.204702518612579e-03 & 1.500937328757639e-03 & 6.561155330342952e-02 & 9.328875093678128e-01\\
\bottomrule
\end{tabular}
\bigskip

\begin{tabular}{rrrrr}
\multicolumn{5}{c}{\rmfamily $d=11$, $n_p=28$, type $[1,3,3]$, quality NI}\\
\toprule
1 & -6.240162943348243e-02 & 3.333333333333333e-01 & 3.333333333333333e-01 & 3.333333333333333e-01\\
3 & ~8.172717083864956e-02 & 4.327128944999808e-01 & 2.836435527500096e-01 & 2.836435527500096e-01\\
3 & ~4.896465250586733e-02 & 7.112097470078995e-01 & 1.443951264960502e-01 & 1.443951264960502e-01\\
3 & ~1.380135388627025e-02 & 9.345387060670882e-01 & 3.273064696645591e-02 & 3.273064696645591e-02\\
6 & ~5.340235179593824e-02 & 1.163964134004934e-01 & 3.353500785250607e-01 & 5.482535080744460e-01\\
6 & ~2.632505873531205e-02 & 2.261151330038821e-02 & 3.724949218910710e-01 & 6.048935648085408e-01\\
6 & ~2.509293909226987e-02 & 2.799022568208098e-02 & 1.649013104719147e-01 & 8.071084638460043e-01\\
\bottomrule
\end{tabular}
\bigskip

\begin{tabular}{rrrrr}
\multicolumn{5}{c}{\rmfamily $d=11$, $n_p=30$, type $[0,2,4]$, quality PI}\\
\toprule
3 & ~5.623165917468111e-02 & 4.470587017120257e-01 & 2.764706491439872e-01 & 2.764706491439872e-01\\
3 & ~4.776140553308587e-02 & 7.160355004191862e-01 & 1.419822497904069e-01 & 1.419822497904069e-01\\
6 & ~5.473425616511274e-02 & 1.230230372259886e-01 & 3.322412718141577e-01 & 5.447356909598537e-01\\
6 & ~2.801405891803870e-02 & 2.440399719079145e-02 & 3.726974607915782e-01 & 6.028985420176303e-01\\
6 & ~2.479873781819909e-02 & 2.785386480846610e-02 & 1.671151632795214e-01 & 8.050309719120125e-01\\
6 & ~7.123081411432649e-03 & 2.715094170950369e-02 & 3.934509469602969e-02 & 9.335039635944666e-01\\
\bottomrule
\end{tabular}
\bigskip

\begin{tabular}{rrrrr}
\multicolumn{5}{c}{\rmfamily $d=11$, $n_p=30$, type $[0,2,4]$, quality PI}\\
\toprule
3 & ~7.217804239720927e-02 & 2.128984074010406e-01 & 3.935507962994797e-01 & 3.935507962994797e-01\\
3 & ~4.321543561536084e-02 & 4.041868382205103e-02 & 4.797906580889745e-01 & 4.797906580889745e-01\\
6 & ~5.817371032216302e-02 & 1.253995635366209e-01 & 2.659762019033016e-01 & 6.086242345600775e-01\\
6 & ~1.697584832298340e-02 & 1.240997015369853e-02 & 2.853641853869646e-01 & 7.022258444593369e-01\\
6 & ~2.759146415644959e-02 & 5.279205798821771e-02 & 1.372353674781709e-01 & 8.099725745336114e-01\\
6 & ~6.228904858785596e-03 & 5.100344564582806e-03 & 5.681715578857245e-02 & 9.380824996468447e-01\\
\bottomrule
\end{tabular}
\bigskip

\begin{tabular}{rrrrr}
\multicolumn{5}{c}{\rmfamily $d=11$, $n_p=30$, type $[0,2,4]$, quality PI}\\
\toprule
3 & ~5.832566212744962e-02 & 4.504543641757660e-01 & 2.747728179121170e-01 & 2.747728179121170e-01\\
3 & ~1.387599563149455e-02 & 9.343162547117883e-01 & 3.284187264410585e-02 & 3.284187264410585e-02\\
6 & ~5.336354239334048e-02 & 1.214249938587573e-01 & 3.347953592492709e-01 & 5.437796468919718e-01\\
6 & ~2.781232090455112e-02 & 2.400046762583091e-02 & 3.718599950903680e-01 & 6.041395372838011e-01\\
6 & ~2.525341262993140e-02 & 1.270006888757827e-01 & 1.580400095523586e-01 & 7.149593015718587e-01\\
6 & ~2.413656185937158e-02 & 2.703971256481997e-02 & 1.649223432616415e-01 & 8.080379441735386e-01\\
\bottomrule
\end{tabular}
\bigskip

\begin{tabular}{rrrrr}
\multicolumn{5}{c}{\rmfamily $d=11$, $n_p=30$, type $[0,2,4]$, quality PI}\\
\toprule
3 & ~4.888313586239229e-02 & 7.118787020151916e-01 & 1.440606489924042e-01 & 1.440606489924042e-01\\
3 & ~1.381436845965494e-02 & 9.345094686995734e-01 & 3.274526565021330e-02 & 3.274526565021330e-02\\
6 & ~3.511195082989262e-02 & 2.404698682473196e-01 & 3.020758710955400e-01 & 4.574542606571404e-01\\
6 & ~5.015309687490970e-02 & 1.100051862098844e-01 & 3.372221080172993e-01 & 5.527727057728163e-01\\
6 & ~2.503977132164279e-02 & 2.152075977107619e-02 & 3.731132768172773e-01 & 6.053659634116466e-01\\
6 & ~2.501309547919794e-02 & 2.787575641695829e-02 & 1.650070701309573e-01 & 8.071171734520844e-01\\
\bottomrule
\end{tabular}
\bigskip

\begin{tabular}{rrrrr}
\multicolumn{5}{c}{\rmfamily $d=12$, $n_p=33$, type $[0,5,3]$, quality PI}\\
\toprule
3 & ~6.254121319590276e-02 & 4.570749859701478e-01 & 2.714625070149261e-01 & 2.714625070149261e-01\\
3 & ~4.991833492806094e-02 & 1.197767026828138e-01 & 4.401116486585931e-01 & 4.401116486585931e-01\\
3 & ~2.426683808145203e-02 & 2.359249810891690e-02 & 4.882037509455416e-01 & 4.882037509455416e-01\\
3 & ~2.848605206887754e-02 & 7.814843446812914e-01 & 1.092578276593543e-01 & 1.092578276593543e-01\\
3 & ~7.931642509973638e-03 & 9.507072731273288e-01 & 2.464636343633559e-02 & 2.464636343633559e-02\\
6 & ~4.322736365941421e-02 & 1.162960196779266e-01 & 2.554542286385173e-01 & 6.282497516835561e-01\\
6 & ~2.178358503860756e-02 & 2.303415635526714e-02 & 2.916556797383410e-01 & 6.853101639063919e-01\\
6 & ~1.508367757651144e-02 & 2.138249025617059e-02 & 1.272797172335894e-01 & 8.513377925102400e-01\\
\bottomrule
\end{tabular}
\bigskip

\begin{tabular}{rrrrr}
\multicolumn{5}{c}{\rmfamily $d=12$, $n_p=33$, type $[0,5,3]$, quality NI}\\
\toprule
3 & ~5.992157930040981e-02 & 4.529711389058645e-01 & 2.735144305470678e-01 & 2.735144305470678e-01\\
3 & ~2.807875643954752e-02 & 2.368917770665134e-02 & 4.881554111466743e-01 & 4.881554111466743e-01\\
3 & ~5.252899601772313e-02 & 7.321646206597614e-01 & 1.339176896701193e-01 & 1.339176896701193e-01\\
3 & ~1.617355627623166e-03 & 9.927521106062486e-01 & 3.623944696875718e-03 & 3.623944696875718e-03\\
3 & -1.062024194350891e-01 & 8.745335893925173e-01 & 6.273320530374134e-02 & 6.273320530374134e-02\\
6 & ~5.491810838782295e-02 & 1.202241316165672e-01 & 3.319346641205961e-01 & 5.478412042628367e-01\\
6 & ~2.544555194057983e-02 & 2.405691547178780e-02 & 2.633690807904016e-01 & 7.125740037378106e-01\\
6 & ~6.833087236315665e-02 & 4.346701716737803e-02 & 7.887216477846390e-02 & 8.776608180541581e-01\\
\bottomrule
\end{tabular}
\bigskip

\begin{tabular}{rrrrr}
\multicolumn{5}{c}{\rmfamily $d=13$, $n_p=37$, type $[1,4,4]$, quality PI}\\
\toprule
1 & ~6.666531183964321e-02 & 3.333333333333333e-01 & 3.333333333333333e-01 & 3.333333333333333e-01\\
3 & ~5.637138317907531e-02 & 1.416807749137785e-01 & 4.291596125431107e-01 & 4.291596125431107e-01\\
3 & ~5.703687953133591e-02 & 5.483198542641960e-01 & 2.258400728679020e-01 & 2.258400728679020e-01\\
3 & ~2.704770288106011e-02 & 2.514270940526753e-02 & 4.874286452973662e-01 & 4.874286452973662e-01\\
3 & ~3.254577710106209e-02 & 7.510843558720015e-01 & 1.244578220639993e-01 & 1.244578220639993e-01\\
6 & ~3.846210380706763e-02 & 7.127447151191104e-02 & 2.845207640198182e-01 & 6.442047644682707e-01\\
6 & ~9.138438814371032e-03 & 4.935323489543055e-03 & 2.862147535443420e-01 & 7.088499229661149e-01\\
6 & ~1.751340205091933e-02 & 2.673280979433629e-02 & 1.245254158513282e-01 & 8.487417743543355e-01\\
6 & ~3.940965341434761e-03 & 1.635078050759145e-02 & 3.285424868085981e-02 & 9.507949708115487e-01\\
\bottomrule
\end{tabular}
\bigskip

\begin{tabular}{rrrrr}
\multicolumn{5}{c}{\rmfamily $d=13$, $n_p=37$, type $[1,4,4]$, quality PI}\\
\toprule
1 & ~6.796003658683164e-02 & 3.333333333333333e-01 & 3.333333333333333e-01 & 3.333333333333333e-01\\
3 & ~5.560196753045333e-02 & 1.461171714803992e-01 & 4.269414142598004e-01 & 4.269414142598004e-01\\
3 & ~5.827848511919998e-02 & 5.572554274163342e-01 & 2.213722862918329e-01 & 2.213722862918329e-01\\
3 & ~2.399440192889473e-02 & 2.184610709492130e-02 & 4.890769464525393e-01 & 4.890769464525393e-01\\
3 & ~6.052337103539172e-03 & 9.569806377823136e-01 & 2.150968110884318e-02 & 2.150968110884318e-02\\
6 & ~3.464127614084837e-02 & 6.801224355420665e-02 & 3.084417608921178e-01 & 6.235459955536756e-01\\
6 & ~2.417903981159382e-02 & 8.789548303219732e-02 & 1.635974010678505e-01 & 7.485071158999522e-01\\
6 & ~9.590681003543263e-03 & 5.126389102382369e-03 & 2.725158177734297e-01 & 7.223577931241880e-01\\
6 & ~1.496540110516567e-02 & 2.437018690109383e-02 & 1.109220428034634e-01 & 8.647077702954428e-01\\
\bottomrule
\end{tabular}
\bigskip

\begin{tabular}{rrrrr}
\multicolumn{5}{c}{\rmfamily $d=13$, $n_p=37$, type $[1,4,4]$, quality NI}\\
\toprule
1 & -1.056360738456401e-01 & 3.333333333333333e-01 & 3.333333333333333e-01 & 3.333333333333333e-01\\
3 & ~9.690034727804083e-02 & 4.242813469723264e-01 & 2.878593265138368e-01 & 2.878593265138368e-01\\
3 & ~5.018216635227328e-02 & 1.131979829755801e-01 & 4.434010085122100e-01 & 4.434010085122100e-01\\
3 & ~2.102855973694512e-02 & 2.147839246570478e-02 & 4.892608037671476e-01 & 4.892608037671476e-01\\
3 & ~2.672599067349385e-02 & 8.000919735991067e-01 & 9.995401320044664e-02 & 9.995401320044664e-02\\
6 & ~4.607891261373841e-02 & 1.192644209390402e-01 & 2.495054689443353e-01 & 6.312301101166245e-01\\
6 & ~2.135290362787687e-02 & 2.361515966854858e-02 & 3.017346932372817e-01 & 6.746501470941697e-01\\
6 & ~1.450394248637925e-02 & 2.089774464177895e-02 & 1.430257663819780e-01 & 8.360764889762430e-01\\
6 & ~4.918388225902278e-03 & 1.497391087216802e-02 & 3.985330690077000e-02 & 9.451727822270620e-01\\
\bottomrule
\end{tabular}
\bigskip

\begin{tabular}{rrrrr}
\multicolumn{5}{c}{\rmfamily $d=13$, $n_p=37$, type $[1,4,4]$, quality NI}\\
\toprule
1 & -9.794282830779282e-01 & 3.333333333333333e-01 & 3.333333333333333e-01 & 3.333333333333333e-01\\
3 & ~3.851843248196896e-01 & 3.731714975085028e-01 & 3.134142512457486e-01 & 3.134142512457486e-01\\
3 & ~4.707074375139372e-02 & 1.036731386947847e-01 & 4.481634306526077e-01 & 4.481634306526077e-01\\
3 & ~1.947159860023236e-02 & 1.876094369614878e-02 & 4.906195281519256e-01 & 4.906195281519256e-01\\
3 & ~8.598559295098043e-03 & 9.482924272447965e-01 & 2.585378637760174e-02 & 2.585378637760174e-02\\
6 & ~4.686871755715881e-02 & 1.280748545621051e-01 & 2.559060504714026e-01 & 6.160190949664922e-01\\
6 & ~2.321458227198339e-02 & 2.466575888342652e-02 & 2.937275637422268e-01 & 6.816066773743466e-01\\
6 & ~1.865886541493881e-02 & 7.660504708654196e-02 & 1.303752945362578e-01 & 7.930196583772002e-01\\
6 & ~1.099993536903349e-02 & 1.517155770635715e-02 & 1.295372633260821e-01 & 8.552911789675607e-01\\
\bottomrule
\end{tabular}
\bigskip

\begin{tabular}{rrrrr}
\multicolumn{5}{c}{\rmfamily $d=13$, $n_p=37$, type $[1,4,4]$, quality NI}\\
\toprule
1 & -6.961938918173517e-01 & 3.333333333333333e-01 & 3.333333333333333e-01 & 3.333333333333333e-01\\
3 & ~2.916465048495467e-01 & 3.802409844193277e-01 & 3.098795077903361e-01 & 3.098795077903361e-01\\
3 & ~4.702284598963952e-02 & 1.047710914472636e-01 & 4.476144542763682e-01 & 4.476144542763682e-01\\
3 & ~2.943707721076129e-02 & 7.979373791902761e-01 & 1.010313104048620e-01 & 1.010313104048620e-01\\
3 & ~7.060478452030929e-03 & 9.526826973365726e-01 & 2.365865133171368e-02 & 2.365865133171368e-02\\
6 & ~4.748213927210297e-02 & 1.255412745019193e-01 & 2.552672945625177e-01 & 6.191914309355630e-01\\
6 & ~1.837479881995625e-02 & 2.035044354178874e-02 & 3.990878803018775e-01 & 5.805616761563338e-01\\
6 & ~1.907290693660829e-02 & 2.679648266402744e-02 & 2.407898961433650e-01 & 7.324136211926076e-01\\
6 & ~1.018568368990188e-02 & 1.643529593541106e-02 & 1.139836801912405e-01 & 8.695810238733484e-01\\
\bottomrule
\end{tabular}
\bigskip
\end{center}

\bibliographystyle{abbrvnat}
\bibliography{Cubature}

\end{document}